\documentclass{amsart}
\usepackage{graphicx}
\newtheorem{thm}{Theorem}[section]

\newtheorem{lem}[thm]{Lemma}
\newtheorem{prop}[thm]{Proposition}
\theoremstyle{definition}
\newtheorem{defn}[thm]{Definition}
\theoremstyle{remark}
\newtheorem{rem}[thm]{Remark}
\numberwithin{equation}{section}

\newcommand{\re}{\rm {Re}}

\def\endproof{\hfill{\vrule height4pt width6pt depth2pt}}
\def\build#1_#2^#3{\mathrel{\mathop{\kern 0pt#1}\limits_{#2}^{#3}}}
\def\tend#1#2{\build\hbox to 12mm{\rightarrowfill}_{#1\rightarrow #2}^{}}
\def\netendpas#1#2{\build\hbox to 12mm{$\not \longrightarrow$}_{#1 \rightarrow
#2}^{}}

\def\tendn{\tend{n} {\infty}}

\def \converge#1#2#3{\build\hbox to 15mm {\rightarrowfill}_{#1\rightarrow #2}^
{\hbox{\scriptsize #3}}}

\newcommand\sym{\fam\comfam\com}
\font\tensym=msbm10 at 12pt \font\sevensym=msbm7
\font\fivesym=msbm5 
\newfam\symfam
\textfont\symfam=\tensym \scriptfont\symfam=\sevensym
\scriptscriptfont\symfam=\fivesym
\renewcommand\sym{\fam\symfam\relax}

\newcommand\N{{\sym N}}
\newcommand\Z{{\sym Z}}
\newcommand\Q{{\sym Q}}

\newcommand\T{{\sym T}}
\newcommand\Prob {{\sym P}}

\def\tower{
\setlength{\unitlength}{1mm}
\begin{picture}(30,30)
\put (0,0) {\framebox (50,20)} \put (3,-4) {\makebox (3,-4)[br]
{$B_{k}$}} \put (16,-6) {\makebox (16,-6)[b]{$\underbrace {
  \quad  \quad   \quad  \quad  \quad  \quad  \quad  \quad \quad  \quad
\quad  \quad  \quad
 }_{B_{k-1}}$}}
\put (17,-16) {\makebox (17,-16)[b]{$\overline {
     \quad  \quad  \quad \quad ~ ~ p_{k-1}
\quad  \quad   \quad  \quad \quad   \quad ~ ~
 }$}}
\multiput (5,0)(5,0){10}{\line (0,20){20}} \put (0,4) {\line
(4,0){50}} \put (3,0){\vector (0,1){4}} \thicklines \put (0,0)
{\line (5,0){5}}
\multiput (0,21)(0,1){5}{\line (5,0){5}} \put (-3,27){\makebox
(-2,27)[bl] {$a_1^{(k-1)}$}} \put (5,21){\line(5,0){5}} \put
(5,22){\line(5,0){5}} \put (8,24){\makebox (8,24)[bl]
{$a_2^{(k-1)}$}} \put (10,21){\makebox (10,21)[b] {$\cdots
\cdots$}} \multiput (25,21)(0,1){6}{\line (25,0){5}} \put
(23,28){\makebox (23,28)[bl] {$a_j^{(k-1)}$}} \put
(26,21){\makebox (26,21)[b] {$\cdots \cdots$}} \multiput
(45,21)(0,1){3}{\line (45,0){5}} \put (45,25){\makebox (45,25)[bl]
{$a_{p_{k-1}}^{(k-1)}$}} \put (-40,-25){\centerline{Figure 1:
$k^{\hbox{th}}$--tower.}}
\end{picture}}

\begin{document}

\title[Ornstein transformations]{A new class of Ornstein transformations with singular spectrum \\
Une nouvelle Classe de transformations d'Ornstein a spectre singulier}%
\author{E. H. EL ABDALAOUI }
\author{F. Parreau }
\author{A. A. PRIKHOD'KO
}

\pagestyle{myheadings} \markboth{El Abdalaoui-Parreau \&
Prikhod'ko} {Generalized Ornstein Transformation}




\dedicatory{\footnotesize {Department of Mathematics, University
of Rouen, LMRS, UMR 60 85, Mont Saint Aignan  76821, France. \\
\footnotesize e-mail : Elhocein.Elabdalaoui@univ-rouen.fr }\\
\footnotesize {Department of Mathematics, University Paris 13
 LAGA, UMR 7539 CNRS ,
 99 Av. J-B Cl\'ement, 93430 Villetaneuse, France}\\
 \footnotesize{email : parreau@math.univ-paris13.fr }\\
\footnotesize{Department of Mathematics, Moscow State University.}\\
\footnotesize{email : apri7@geocities.com}
 }
\begin{abstract}
It is shown that for any family of probability measures in
Ornstein type constructions the corresponding transformation has
almost surely a singular spectrum. This is a new generalization of
Bourgain's theorem \cite{Bourgain}, the same result is proved for
Rudolph's construction \cite{Rudolph}.\\

\noindent{}{R\'ESUM\'E.} On montre que pour toute famille de
mesures de probabilit\'es dans la construction
d'Ornstein, les transformations r\'esultantes ont un spectre
presque s\^urement singulier. On obtient ainsi une nouvelle
g\'en\'eralisation d'un th\'eor\'eme d\^u \`a Bourgain
\cite{Bourgain}. Un r\'esultat similaire est obtenu pour les
transformations de Rudolph \cite{Rudolph}.
\end{abstract}
\maketitle
\section{Introduction}
In this note we investigate the spectral analysis of a generalized
class of Ornstein transformations. There are several
generalizations of Ornstein transformations. Here we are concerned
with arbitrary product probability space associated to random
construction of the family of rank one transformations.
Namely, in the Ornstein's construction, the probability
space is equipped with the infinite product of uniform
probability measures on some finite subsets of
$\mathbb{Z}$. Here, the probability space is
equipped with the infinite product of probability measures
$(\xi_m)_{m \in \N}$ on a family $(X_m)_{m \in \N}$ of finite
subsets of $\mathbb{Z}$. We establish that for any choice of the
family $(\xi_m)_{m \in \N}$ the associated Ornstein
transformations has almost surely singular spectrum.\\

Let us recall that Ornstein introduced these transformations in
1967 in \cite{Ornstein} and proved that the
mixing property occurs almost surely. Until 1991, these
transformations which have simple spectrum appeared as a
candidate for an affirmative answer to Banach's well-known
problem whether a dynamical system $(X,\mathcal{B},\mu)$ may have
simple Lebesgue spectrum. But, in 1991, J. Bourgain in
\cite{Bourgain}, using Riesz products techniques, proved that
Ornstein transformations have almost surely singular
spectrum. Subsequently, I. Klemes \cite{Klemes}, I. Klemes \& K. Reinhold
\cite{Klemes-Reinhold} obtain that the spectrum of the mixing
subclass of staircase transformations of T. Adams \cite{Adams} and
T. Adams \& N. Friedman \cite{Adams2} have singular spectrum. They
conjectured that rank one transformations allways have singular spectrum.\\

In this paper, using the techniques of J. Bourgain generalized in
\cite{elAbdal1}, we extend Bourgain's theorem to the generalized
Ornstein transformations associated to a large family of random
constructions.\\

Firstly, we shall recall some basic facts from spectral
theory. A nice account can be found in the appendix of
\cite{Parry}. We shall assume that the reader is
familiar with the method of cutting and stacking for constructing
rank one transformations.\\


Given $T:(X,\mathcal{B},\mu)\mapsto (X,\mathcal{B},\mu)$ a measure preserving invertible
transformation and denoting by $U_T f$ the operator  $U_T
f(x)=f(T^{-1}x)$ on $L^2(X,\mathcal{B},\mu))$, recall that to any $f\in L^2(X)$ there
corresponds a positive measure $\sigma_f$ on $\T$, the
unit circle, defined by $\hat{\sigma}_f(n)= <U_T^nf,f>$. \\

\begin{defn}
The maximal spectral type of $T$ is the equivalence class of Borel
measures $\sigma$ on $\T$ (under the equivalence relation $\mu_1
\sim \mu_2$ if  and only if $\mu_1<<\mu_2$ and $\mu_2<<\mu_1$),
such that
 $\sigma_f<<\sigma$ for all $f\in L^2(X)$ and
if $\nu$ is another measure for which $\sigma_f<<\nu$
for all $f\in L^2(X)$ then $\sigma << \nu$.\\
\end{defn}

There exists a Borel measure $\sigma=\sigma_f$ for some $f\in L^2(X)$, such that $\sigma$ is in
the equivalence class defining the maximal spectral type of $T$.
By abuse of notation, we will call this measure the maximal
spectral type measure. The reduced maximal type $\sigma_0$ is the
maximal spectral type of $U_T$ on $L_0^2(X)\stackrel{\rm
{def}}{=}\{f \in L^2(X)~:~ \displaystyle \int f d\mu=0 \}$. The
spectrum of $T$ is said to be discrete (resp. continuous, resp.
singular, resp. absolutely continuous , resp. Lebesgue ) if
$\sigma_0$ is discrete ( resp. continuous, resp. singular, resp.
absolutely continuous with respect to the Lebesgue measure or
equivalent to the Lebesgue measure).  We write
\[Z(h) \stackrel{\rm {def}}{=} \overline {{\rm {span}} \{U_T^nh,n \in \Z \} }.
\]
\noindent{}$T$ is said to have simple spectrum, if there exists $h
\in L^2(X)$ such that
\[Z(h) =L^2(X).\]

\section{Rank One Transformation by Construction}

\noindent Using the cutting and stacking method described in
\cite{Friedman1}, \cite{Friedman2}, one defines inductively a
family of measure preserving transformations, called rank one
transformations, as follows
\vskip 0.1cm
Let $B_0$ be the unit
interval equipped with the Lebesgue measure. At stage one we
divide $B_0$ into $p_0$ equal parts, add spacers and form a stack
of height $h_{1}$ in the usual fashion. At the $k^{th}$ stage we
divide the stack obtained at the $(k-1)^{th}$ stage into $p_{k-1}$
equal columns, add spacers and obtain a new stack of height
$h_{k}$. If during the $k^{th}$ stage of our construction  the
number of spacers put above the $j^{th}$ column of the
$(k-1)^{th}$ stack is $a^{(k-1)}_{j}$, $ 0 \leq a^{(k-1)}_{j} <
\infty$,  $1\leq j \leq p_{k}$, then we have
\begin{eqnarray*}
h_{k} &=& p_{k-1}h_{k-1} +  \sum_{j=1}^{p_{k-1}}a_{j}^{(k-1)}, ~~~\forall k \geq 1,\\
h_0&=&1.
\end{eqnarray*}
\vskip 1.5 cm \hskip 3.5cm \tower \vskip 3.0cm
\noindent{}Proceeding in this way we get a rank one transformation
$T$ on a certain measure space $(X, \mathcal{B} ,\nu)$ which may
be finite or
$\sigma-$finite depending on the number of spacers added. \\
\noindent{} The construction of any rank one transformation thus
needs two parameters $(p_k)_{k=0}^\infty$ (cutting parameter)
and $((a_j^{(k)})_{j=1}^{p_k})_{k=0}^\infty$
(spacers parameter). We put

$$T \stackrel {def}= T_{(p_k, (a_j^{(k)})_{j=1}^{p_k})_{k=0}^\infty}$$

\noindent In \cite{Bourgain},\cite{Choksi-Nadkarni} and
\cite{Klemes-Reinhold} it is proved that up to some discret
measure, the spectral type of this transformation is given by
\begin{eqnarray}\label{eqn:type1}
d\sigma  =W^{*}\lim \prod_{k=1}^n\left| P_k\right| ^2d\lambda .
\end{eqnarray}
\begin{eqnarray*}
&&~~{\rm {~where~~}}P_k(z)=\frac 1{\sqrt{p_k}}\left(
\sum_{j=0}^{p_k-1}z^{-(jh_k+\sum_{i=1}^ja_i^{(k)})}\right)\nonumber  \\
\nonumber
&&
\lambda \ {\rm {~denotes~the~normalized~Lebesgue~measure~on~torus~}{\fam%
\symfam\relax T}}{\rm {\ .}} \nonumber \\
\nonumber
&& W^{*}{\rm {~denotes~weak~convergence~on~the~space~of~bounded~Borel~measures~on~}}{\fam%
\symfam\relax T}. \nonumber\\
&&
\nonumber
\end{eqnarray*}
\noindent The polynomials $P_k$ appear naturally from the
induction relation between the bases $B_k$. Indeed
\[
B_k=B_{k+1}\cup T^{h_k+s_k(1)}B_{k+1}\cup \ldots \cup
T^{(p_k-1)h_k+s_k(p_k-1)}B_{k+1},
\]

\[
\nu (B_k)=p_k\nu (B_{k+1}),
\]

\noindent where $s_k(n)=a_1^{(k)}+\ldots +a_n^{(k)}$ and
$s_k(0)=0.$

\noindent Put
\[
f_k=\frac 1{\sqrt{\nu (B_k)}} \chi_{B_k},
\]

\noindent that is the indicator function of the $k$th-base
normalized in the $L^{2}$-norm. So
\[
f_k=\ P_k(U_T)f_{k+1},
\]

\noindent{}where $U_T~:~L^2(X) \longrightarrow L^2(X)$ is defined by $%
\displaystyle U_T(f)(x)=f(T^{-1}x)$. Iterating this relation, we
have
\[
d\sigma _k=\left| P_k\right| ^2d\sigma _{k+1}=\ldots
=\prod_{j=0}^{m-1}\left| P_{k+j}\right| ^2d\sigma _{k+m},
\]

\noindent Where $\sigma_p$ is the spectral measure of $f_p$, $p \geq 0$. \\

\section{ Generalized Ornstein's Class of Transformations}

In Ornstein's construction, the $p_k$'s are rapidly increasing,
and the number of spacers, $a_i^{(k)}$, $1 \leq i\leq p_k-1$, are
chosen randomly. This
may be organized in differently ways as pointed by J. Bourgain in
\cite{Bourgain}. Here, We suppose given  $(t_k)$, $(p_k)$ a sequences of positive
integers and ($\xi_k$) a sequence of probability measure such that the support of each
$\xi_k$ is a subset of $X_k =
\{-\displaystyle \frac{t_k}{2},\cdots,\displaystyle
\frac{t_k}{2}\}$. We choose now independently, according to $\xi_k$
the numbers $(x_{k,i})_{i=1}^{p_k-1}$, and $x_{k,p_k}$ is chosen
deterministically in $\N$. We put, for $1 \leq i \leq p_k$,

$$a_i^{(k)} = t_{k} + x_{k,i} - x_{k,i-1}, ~~{\rm with} ~~x_{k,0}
= 0.$$\\

\noindent{} It follows

$$h_{k+1} = p_k(h_k + t_{k}) + x_{k,p_k}.$$

\noindent{}So the deterministic sequences of positive integers
$(p_k)_{k=0}^\infty$, $(t_k)_{k=0}^\infty $ and
$(x_{k,p_k})_{k=0}^\infty$ determine completely the sequence of
heights $(h_k)_{k=0}^\infty$. The total measure of the resulting
measure space is finite if
\begin{eqnarray}\label{eqn:fini}
\sum_{k=0}^{\infty}\frac{t_k}{h_k}+\sum_{k=0}^\infty
\frac{x_{k,p_k}}{p_kh_k} < \infty.
\end{eqnarray}
 \noindent{}We will assume that this
requirement is satisfied.\\
We thus have a probability space of Ornstein transformations
$\Omega=\prod_{l=0}^\infty X_l^{p_l-1}$ equipped with the natural
probability measure $\Prob \stackrel {\rm def}
{=}\otimes_{l=1}^{\infty} P_l$, where $P_l\stackrel {\rm def}
{=}\otimes_{j=1}^{p_l-1}{\xi_l}$; ${\xi_l}$ is the probability
measure on $X_l$. We denote this space by $(\Omega, {\mathcal
{A}}, {\Prob})$. So $x_{k,i}$, $1 \leq i \leq p_k -1$, is the
projection from $\Omega$ onto the $i^{th}$ co-ordinate space of
$\Omega_k \stackrel {\rm def} {=} X_k^{p_k-1}$, $1 \leq i \leq
p_k-1$. Naturally each point $\omega =(\omega_k =
(x_{k,i}(\omega))_{i=1}^{p_k-1})_{k=0}^\infty$ in $\Omega$ defines
the spacers and therefore a rank one transformation
$T_{\omega,x}$, where $x=(x_{k,p_k})$.

\noindent{}The definition above gives a more general definition of
random construction due to Ornstein. In the particular case of
Ornstein's transformations constructed in \cite{Ornstein},
$t_k=h_{k-1}$, $\xi_k$ is uniform distribution and  $p_k >> h_{k-1}$. \\
\noindent{} We recall that Ornstein in \cite{Ornstein} proved that
there exist a sequence ${(p_k,x_{k,p_k})}_{k \in \N}$ such that,
$T_{\omega,x}$ is almost surely mixing. Later in \cite{Prikhodko},
Prikhod'ko obtain the same result for some special choice of the
sequence of the distribution ${(\xi_m)}$ and recently, using the
idea of D. Creutz and C. E. Silva \cite{Creutz-Silva} one can
extend this result to a large class of the family of the
probability measure associated to Ornstein construction. In our
general construction, according to (\ref{eqn:type1}) the spectral
type of each $T_{\omega}$ , up to a discrete measure, is given by
  \[
\sigma_{T_\omega }=\sigma^{(\omega)}
_{\chi_{B_0}}=\sigma^{(\omega)} =W ^{*}\lim \prod_{l=1}^N\frac
1{p_l}\left| \sum_{p=0}^{p_l-1}z^{p(h_l+t_l)+x_{l,p}}\right|
^2d\lambda.
\]

\noindent{} With the above notation, we state our main result\\

\bigskip
\begin{thm}
For every choice of $(p_k), (t_k), (x_{k,p_k})$ and for any family
of probability measures ${\xi_m}$ on finite subset $X_m$ of $\Z$,
$ {m \in \N^*}$. The associated generalized Ornstein
transformations has almost surely singular spectrum. i.e.
\[
\Prob\{\omega~:~ \sigma^{(\omega)} ~\bot ~\lambda\}=1.
\]
\noindent{}Where $\Prob \stackrel {\rm def}
{=}\otimes_{l=0}^{\infty} \otimes_{j=1}^{p_l-1}{\xi_l}$; is the
probability measure on $\Omega=\prod_{l=0}^\infty X_l^{p_l-1}$,
$X_l$ is finite subset of $\Z$.
\end{thm}

\noindent{}Before proceeding to the proof, we remark that it is an
easy exercise to see that the spectrum of Ornstein's
transformation is always singular if the cutting parameter $p_k$
is bounded. In fact, Klemes-Reinhold proved moreover that if
$\displaystyle \sum_{k=0}^{\infty} \frac{1}{{p_k}^2}=\infty$ then
the associated rank one transformation is singular. Henceforth, we
assume that the series $\displaystyle \sum_{k=0}^{\infty} \frac{1}{{p_k}^2}$ converges.\\

\noindent{}We shall adapt Bourgain's proof. For that, we need a
local version of the singularity criterion used by Bourgain. Let
$F$ be a Borel set then with the above notations, we will state
local singularity criterion in the following form\\

\begin{thm}{{\bf (Local Singularity Criterion (LSC))}}
The following are equivalent\\

\item[(i)]  {\it $\sigma_F  \perp \lambda ,$ }, where $\sigma_F=\chi_F.d\sigma$, $\chi_F$ is a indicator
function of $F$.\\

\item[(ii)] $\displaystyle \int_F \prod_{l=1}^n\left|{P_{l}(z)}\right| d\lambda \tendn 0.$\\

\item[(iii)]{\it $\inf \{\displaystyle
\int_F \prod_{l=1}^k\left|{P_{n_l}(z)}\right| d\lambda, ~k\in
{\fam\symfam\relax N},~n_1<n_2<\ldots <n_k\}=0.$ }
\end{thm}

\noindent{}One can adapt the proof of theorem 4.3 in
\cite{Klemes-Reinhold}, or in \cite{elAbdal1}, \cite{Nadkarni}, in
the more general setting.\\
Now, using Lebesgue's dominated convergence theorem and the LSC
, we obtain  \\

\begin{prop}
 The following are equivalent

\item[(i)]  {\it ${\sigma_F ^{(\omega )}}\perp \lambda$ $%
\qquad {\fam\symfam\relax P}~a.s.$ }

\item[(ii)] $\displaystyle \int_F \prod_{l=1}^n\left|{P_{l}(z)}\right| d\lambda d{\fam\symfam\relax P}
\tendn 0.$\\

\item[(iii)] {\it $\inf \{\displaystyle \displaystyle \int \displaystyle %
\int_F \prod_{l=1}^k\left|{P_{n_l}(z)}\right| d\lambda d{\fam\symfam\relax P},~k\in {\fam\symfam\relax N}%
,~n_1<n_2<\ldots <n_k\}=0.$ }
\end{prop}

\noindent{}Fix some subsequence ${\mathcal N=}\left\{
n_1<n_2<\ldots <n_k\right\}$ ,\linebreak[0]$\ k \in
{\fam\symfam\relax N}$, $m > n_k$ and put
\begin{eqnarray*}
Q\left(z\right) &=& \prod_{i=1}^k\vert{P_{n_i}(z)}\vert.
\end{eqnarray*}
\noindent{}Following \cite{Bourgain} ( see also \cite{Klemes}
or in the more general setting \cite{elAbdal3}), we have.\\
\begin{lem}
\[
\displaystyle \!\!\int_F Q \left| P_m\right|~ d\lambda ~ \leq
\frac 12\left( \displaystyle \!\!\int_F  Q d \lambda+\displaystyle
\int_F Q \left| P_m\right|^2 d\lambda\right) -\frac 18\left(
\displaystyle \int_F Q \left| \left| P_m\right| ^2-1 \right|
d\lambda \right)^2.
\]\\
\end{lem}

\noindent{}Now, we assume that $F$ is closed set, it follows

\begin{lem}
$\displaystyle \limsup_{m \rightarrow
\infty}\displaystyle \int_F \!\!\!\! Q \left| {P_m(z)}%
\right| ^2d\lambda(z)  \leq \displaystyle \displaystyle \int_F Q $
$d\lambda(z).$ \\
\end{lem}

\begin{proof}
\noindent{}Observe that the
sequence of probability measures $\left| {P_m(z)}%
\right|^2d\lambda(z)$ converges weakly to the Lebsegue measure.
Then the lemma follows from the classical portmanteau theorem
\footnote{see for example
  \cite{Dudley}. We note that the space $\Omega$ is equiped with the
  standard product topology.} and the proof is complete.
  \end{proof}
  \\

\noindent{}From the lemmas 3.4 and 3.5 we get the following

\begin{lem}
\[
\liminf\displaystyle \int \!\!\!\!\displaystyle \int_F Q \left |
{P_m} \right | d\lambda d{\fam\symfam\relax P} \leq %
\displaystyle \int \!\!\!\!\displaystyle \int_F Q d\lambda d{\fam\symfam\relax P}-
\frac 18\left( {\limsup }%
\displaystyle \int \!\!\!\!\displaystyle \int_F Q \left| \left| {P_m}%
\right| ^2-1\right|  d\lambda d{\fam\symfam\relax P}\right) ^2.
\]
\end{lem}

\noindent{}Clearly,  we need to estimate the quantity
\begin{eqnarray}
\displaystyle \int \displaystyle \int_F  Q \left| \left|{P_m(z)}%
\right| ^2-1\right| d\lambda(z) d{\fam\symfam\relax P}.
\end{eqnarray}
\noindent{}For that, following Bourgain we shall prove the
following

\begin{prop} There exists an absolute constant $K>0$ such that
\begin{eqnarray*}
\limsup
\displaystyle \int \displaystyle \int_F Q \left| \left| {P_m}%
\right| ^2-1\right| d\lambda d\Prob  \geq K  {(\int \displaystyle
\int_F Q
 d\lambda d\Prob -{\liminf }\displaystyle \int \displaystyle
 \int_F Q(z)\phi_m(z) d\lambda d\Prob)}^2,
 \end{eqnarray*}
 \noindent{} where $ \displaystyle \phi_m(z)=
\left |\sum_{p=-\frac{t_m}2}^{p=\frac{t_m}2}
\xi_m(p) z^p \right|^2,$  $z \in \T$\\
\end{prop}
\noindent{}We shall give the proof of proposition 3.7 in the
following section.\\

\section{Khintchine-Bonami inequality}

\noindent Fix $z \in \T$ and $m \in \N^{*}$. Define $\tau $ and
$(\tau_p)_{p=1}^{p_m-1}$ by :
\[
\begin{array}{c}
{\rm ~} \quad \quad \quad \tau :{\fam\symfam\relax Z}\longrightarrow {\fam%
\symfam\relax T} \\
\qquad {\ \qquad \qquad }s\longmapsto z^s.{\rm ~}
\end{array}
\]
\noindent $\tau _p$ is given by $\tau _p=\tau \circ x_{m,p},$
$x_{m,p}$ is the $p^{\mbox {th}}$ projection on $\Omega
_m=X_m^{p_m-1}. $ So
\[
\left| P_m(z)\right| ^2-1={\sum }_{p\neq q}a_{pq}~\tau _p(\omega )
\overline{\tau _q}(\omega ).{\rm {\quad }}~ {\rm {where }}~~\displaystyle %
a_{pq}=\frac{z^{(p-q)(h_m+t_{m})}}{p_m},
\]

\noindent The random variables $(\tau _p)_{p=1}^{p_m-1}$ are
independent. Put
\begin{equation}
\tau _p^{\circ }=\tau _p-\displaystyle \int \tau_p \
d{\fam\symfam\relax P}, ~~~p=1,\cdots,p_m-1.
\end{equation}
and write
\begin{equation}
\begin{array}{c}
\displaystyle \sum a_{pq}~\tau _p~\overline{\tau _q}= \\
\left ( \displaystyle \sum a_{pq} \right )\left| \displaystyle \int \tau_1
\right|^2+\sum
a_{pq} \left (\displaystyle \left(\int \overline{\tau_1 }\right)~\tau _p^{\circ }+\displaystyle %
\left (\int \tau_1 \right )~\overline{\tau_q^{\circ }}\right )+\sum a_{pq}\tau _p^{\circ } ~%
\overline{\tau _q^{\circ }}.
\end{array}
\end{equation}
\noindent Now, using the same arguments as J. Bourgain, let us
consider a random sign $\varepsilon =\left\{ \varepsilon _1,\ldots
,\varepsilon _{p_m-1}\right\} \in \left\{ -1,1\right\} ^{p_m-1}$,
and the probability space
\[
Z_m=\Omega _m\times \left\{ -1,1\right\} ^{p_m-1},~{\rm {where}}~~
\Omega _m=\left\{ -\frac{t_{m}}2,...,\frac{t_{m}}2\right\}
^{p_m-1}.
\]

\noindent Taking the conditional expectation of the following
quantity
\[
\sum a_{pq}\left( \displaystyle \int \tau _p^{\circ}~~ \overline{\tau_1 }+%
\displaystyle \int \overline{\tau _q^{\circ }}~~\tau_1 ~~\right)
+\sum a_{pq}\tau _p^{\circ }~~{\ }\overline{\tau _q^{\circ }}
\]
\noindent with respect to the $\sigma-$algebra ${\mathcal
{B}_{\varepsilon}}$ given by the cylindres sets $A(I,x)$ where
$I\subset \left\{ 1,\ldots ,p_m-1\right\} $, $x\in \Omega _m$ and
\[
A(I,x)=\prod_{i\in I}\left\{ x_i\right\} \times \left\{ -\frac{t_{m}}2,...,%
\frac{t_{m}}2\right\} ^{\left| I^c\right| }\times \left\{
1\right\} ^{\left| I\right| }\times \left\{ -1\right\} ^{\left|
I^c\right| }.
\]

\noindent ($I$ corresponds to $\displaystyle\varepsilon _i=1$,
$\forall i\in I$ and $\varepsilon _i=-1$, $\forall i\notin I$).
\noindent In other words, taking conditional expectation with
respect to the random variables $\tau_p $ for which
$\varepsilon_p=1,$ one finds the following polynomial expression
in $\varepsilon $ of degree 2
\begin{equation}
\sum a_{pq}\left( \frac{1+\varepsilon _p}2\int \overline{\tau_1
}~{\ }\tau _p^{\circ }+\frac{1+\varepsilon _q}2 \int \tau_1 ~{\
}\overline{\tau _q^{\circ }}\right) +\sum
a_{pq}\frac{1+\varepsilon _p}2\frac{1+\varepsilon _q}2\tau
_p^{\circ }~{\ }\overline{\tau _q^{\circ }}
\end{equation}
\noindent So
\begin{equation}
\begin{array}{c}
\displaystyle \int \left| \left|P_m(z)\right| ^2-1|
\right|d{\fam\symfam\relax P} =\displaystyle \int \displaystyle
\int {\fam\symfam\relax E}(\left| \left|
P_m(z)\right| ^2-1 \right|_{|{\mathcal {B}}_{\varepsilon}})~d{\fam%
\symfam\relax P} d\varepsilon\\
\geq \displaystyle \int \displaystyle \int \left| {\fam\symfam\relax E}%
(\left| P_m(z)\right|^2-1 _{|{\mathcal {B}}_{\varepsilon}%
})~~\right| d{\fam\symfam\relax P} d\varepsilon.
\end{array}
\end{equation}

\noindent It follows, by the Khintchine-Bonami inequality, \footnote{%
One can extend easily this inequality to bounded sequences of
independent real random variables, with vanishing expectation.}
\cite{Bonami}, that there exists a positive constant $K$ such that
\begin{equation}
\begin{array}{c}
\displaystyle \int \displaystyle \int \left| {\fam\symfam\relax
E}(\left|
P_m(z)\right| ^2-1 _{|{\mathcal {B}}_{\varepsilon}}]~~\right| d{%
\fam\symfam\relax P} d\varepsilon \\
\geq K \displaystyle \int (\displaystyle \int \left| {\fam\symfam\relax E}%
(\left| P_m(z)\right| ^2-1 _{|{\mathcal {B}}_{\varepsilon}%
}]~~\right|^2 d\varepsilon)^{\frac{1}{2}} d{\fam\symfam\relax P} \\
=K\displaystyle \int {\displaystyle \left (\sum_{p \neq q} \left
|a_{pq}(z){\tau_p^{\circ}(z)} {\overline {\tau_q^{\circ}(z)}}
\right|^2 \right)^{\frac{1}{2}}d{\fam\symfam\relax P}}.
\end{array}
\end{equation}
\noindent{}But all these random variables are bounded by 2. Hence
\begin{equation}
\begin{array}{c}
\displaystyle \int \left| \left| P_m(z)\right| ^2-1\right| d {\fam\symfam\relax P} \\
\geq K'\displaystyle \int
{\displaystyle\frac{1}{p_m^2}}\displaystyle\sum \left|{\tau_p^{\circ}}(z)
{\overline
{\tau_q^{\circ}(z)}}\right |^2 d{\fam\symfam\relax P} \\
=K^{\prime}{\displaystyle\frac{1}{p_m^2}}\displaystyle\sum \left (\displaystyle \int
\left |\tau_p^{\circ}(z)\right|^2d{%
\fam\symfam\relax P}\right)^2 \\
= K^{\prime}{\displaystyle\frac {(p_m-1)(p_m-2)}{{p_m}^2}}\left(\displaystyle \int \left|\tau_1^{\circ}(z)
\right|^2d%
{\fam\symfam\relax P}\right)^2.
\end{array}
\end{equation}

\noindent{}Since
\begin{equation}
\begin{array}{c}
\displaystyle \int |\tau_1^{\circ}(z)|^2d{\fam\symfam\relax
P}=var(\tau_1(z))
\\
=1-
\left|\displaystyle\sum_{s=-\frac{t_m}2}^{\frac{t_m}2}\xi_m(s)z^s\right|
^2.
\end{array}
\end{equation}
Now, combined (4.6) with (4.7) to obtain\\
\begin{equation}
\begin{array}{c}
\displaystyle \int \left| \left| P_m(z)\right| ^2-1\right| d {\fam\symfam\relax P} \\
\geq K'\displaystyle \frac{(p_m-1)(p_m-2)}{{{p_m}^2}}
\displaystyle  { \left (1-
\left|\displaystyle\sum_{s=-\frac{t_m}2}^{\frac{t_m}2}\xi_m(s)z^s\right|
^2\displaystyle  \right )}^2
\end{array}
\end{equation}
Finally, Multiply (4.8) by
\begin{equation}
\displaystyle \int \prod_{j\in {\mathcal {N}}}\left|P_j(z)\right| d{\fam%
\symfam\relax P}.
\end{equation}
\noindent{}Using the independence of (4.9) and $\displaystyle %
|1-|P_m(z)|^2|$. Integrating over $F$ with respect to the Lebesgue
measure to get
\begin{eqnarray}{\label{eq:eq1}}
\nonumber\int_{\Omega} \int_F Q \left| \left| {P_m(z)}%
\right| ^2-1\right| d\lambda d\Prob \\\geq  K' (\int_{\Omega}
\int_F Q (1-\phi_m(z))^{2} d\lambda d\Prob)
\end{eqnarray}
\noindent{}where
$\phi_m(z)=\left|\displaystyle\sum_{s=-\frac{t_m}2}^{\frac{t_m}2}\xi_m(s)z^s\right|^{2
} $. Apply Cauchy-Schwarz inequality to obtain
\begin{eqnarray}{\label{eq:eq2}}
\nonumber \int_{\Omega} \int_F Q (1-\phi_m(z)) d\lambda d\Prob &
\leq &{\left ( \int_{\Omega} \int_F Q d\lambda d\Prob
\right)}^{\frac12} {\left (\int_{\Omega} \int_F Q
(1-\phi_m(z))^{2} d\lambda
d\Prob\right)}^{\frac12}\\
& \leq & {\left (\int_{\Omega} \int_F Q (1-\phi_m(z))^{2} d\lambda
d\Prob\right)}^{\frac12}
\end{eqnarray}
\noindent{}Combined (\ref{eq:eq1}) and (\ref{eq:eq2}) and take
liminf to finish the proof of the proposition 3.7.
\endproof\\



\noindent{}Now, passing to a subsequence we may assume that
$\phi_m$ converge weakly in $L^{2}(\lambda)$ to some function
$\phi$ in $L^2(\lambda)$. Then,
\[
\widehat{\phi}(n)=\lim_{m\longrightarrow
{\infty}}\widehat{\phi_m}(n)\geq 0.~~~~~{\rm {for~any~}} n \in \Z,
\]
\noindent{}and
\[
\sum_{n}\widehat{\phi}(n) \leq 1.
\]
\noindent{}Hence, the Fourier series of $\phi$ converge absolutely
and we may assume
\[
\phi(z)=\sum_{n}\widehat{\phi}(n)z^n,
\]
\noindent{}In particular $\phi$ is a continuous function. We
deduce that the set $\{\displaystyle \phi(z)=1\}$ is either the
torus or a finite subgroup of the torus.\\

\begin{rem}
It is any easy exercise to see that if the set $\{\displaystyle
\phi=1\}$ is not a null set with respect to Lebesgue measure then,
for any $z \in \T$,
\begin{eqnarray*}
\phi(z)=1\\
\textrm{~and~}
 \max_{s \in X_{m}} \xi_{m}(s) \tend{m}{\infty} 1.
 \end{eqnarray*}
\end{rem}

 \noindent{}We shall, now, prove, our main result in the
following sections.

\section{ On the Ornstein probability space for which  $\underline{\lim }\max_{s \in
  X_{m}}\xi_{m}(s)<1$}

In this section, we assume that $\underline{\lim }\max_{s \in
  X_{m}}\xi_{m}(s)<1$. So, we may choose $\phi$ the weak limite of
  subsequence of $\phi_m$ so that $\widehat{\phi}(0)<1$ and
  $\{\phi=1\}$ is a finite. Let $\varepsilon>0$, put
\[
F_{\varepsilon}\stackrel{def}{= }\{z \in \T ~:~ 1-\phi(z) \geq
\varepsilon\}.
\]
 \noindent{}We get easily that $F_{\varepsilon}$ is a closed set and
 we have also the following proposition\\
\begin{prop} There exists an absolute constant $K>0$ such that
$$
\underline{\lim }%
\displaystyle \iint_{F_{\varepsilon}} Q \left| \left| {P_m(z)}%
\right| ^2-1\right| d\lambda \Prob \geq K{\varepsilon}^{2} {\left
( \iint_{F_{\varepsilon}} Q d\lambda \Prob\right)}^{2}.
$$\\
\end{prop}

\begin{proof}
\noindent{}Apply the proposition 3.7 to get that there exists an
constant $K>0$ for which we have
\begin{eqnarray*}
&&\liminf \displaystyle \iint_{F_{\varepsilon}} Q   \left |\left|
P_m(z) \right| ^2-1\right| d \lambda \Prob \\
&\geq& K {\left
(\iint_F Q d\lambda d\Prob -\overline{\lim }\displaystyle \iint_F
Q(z) \left |\sum_{p=-\frac{t_m}2}^{p=\frac{t_m}2} \xi_m(p) z^p \right |^2
\lambda
d\Prob \right)}^{2}\\
&\geq& K {\left (\iint_{F_{\varepsilon}} Q \left (1-\phi(z) \right) \lambda
d\Prob \right)}^{2}\\
&\geq& K \varepsilon^{2} {\left ( \iint_{F_{\varepsilon}} Q
\lambda d\Prob \right)}^{2}.
\end{eqnarray*}
 \noindent{} The proof of the proposition is complete.
\end{proof}\\

\noindent {\bf Proof of the theorem 3.1.in the case of
$\underline{\lim }\max_{s \in
  X_{m}}\xi_{m}(s)<1$}\\

\noindent{}First, for fixed $\varepsilon>0$, let us choose the
good subsequence ${\mathcal {N}}\stackrel {def}{=}\{n_k, k \geq 0
\}$. Observe that from the propositions 3.6. and 5.1. one can
write
\[
\overline{\lim }\displaystyle \int \int_{F_{\varepsilon}} Q
 \left | {P_m(z)}\right | d\lambda(z) d{\fam\symfam\relax
P} \leq
\displaystyle \int \int_{F_{\varepsilon}} Q -\frac 18K^2\varepsilon^4 \left(\displaystyle \int %
\displaystyle \int_{F_{\varepsilon}} Q d\lambda
d{\fam\symfam\relax P}\right)^4,
\]
\noindent and from this last inequality we shall construct
${\mathcal {N}}$. In
fact, suppose we have chosen the $k$ first elements of the subsequence $%
{\mathcal {N}}$. We wish to define the ${(k+1)}^{{\rm {th}}}$
element. Let $m>n_k $ such that
\[
\displaystyle \int \displaystyle \int_{F_{\varepsilon}}   Q
 \left | {P_m(z)} \right | d\lambda(z) d{\fam\symfam\relax
P} \leq
\displaystyle \int \displaystyle \int_{F_{\varepsilon}}  Q d\lambda d{%
\fam\symfam\relax P} -\frac 18K^2 \varepsilon^4
\left(\displaystyle \int \displaystyle \int_{F_{\varepsilon}}  Q
d\lambda d{\fam\symfam\relax P}\right)^4,
\]
\noindent{}and put $n_{k+1}{=}m.$ It follows that the elements of
the subsequence ${\mathcal {N}}$ verify
\[
\begin{array}{c}
\displaystyle \displaystyle \int \displaystyle
\int_{F_{\varepsilon}} \prod_{i=1}^{k+1}|{P_{n_i}(z)}| d\lambda
d{\fam\symfam\relax P} \leq \displaystyle \displaystyle \int
\displaystyle \int_{F_{\varepsilon}} \prod_{i=1}^k|{P_{n_i}(z)}|
d\lambda d{\fam\symfam\relax P}-
\displaystyle \frac18 K^2 \varepsilon^4{\left(\displaystyle \displaystyle \int \displaystyle %
\int_{F_{\varepsilon}} \prod_{i=1}^k|{P_{n_i}(z)}| d\lambda d
{\fam\symfam\relax P}\right)^4}.
\end{array}
\]
\noindent{}We deduce that the sequence $(\displaystyle \int
\displaystyle \int_{F_{\varepsilon}} \prod_{i=1}^k|{P_{n_i}(z)}|
d\lambda {\fam\symfam\relax P})_{k \geq 1}$ is decreasing and
converges to the limit $l_{\varepsilon}$ which verifies
\[
l_{\varepsilon}\leq l_{\varepsilon}-\frac18
K^2\varepsilon^4l_{\varepsilon}^4,
\]
\noindent{}and this implies that $l_{\varepsilon}=0$. Hence,
${\sigma_{F_{\varepsilon}} ^{(\omega )}}$ is singular. But,
\[ \bigsqcup_{\varepsilon>0, \varepsilon \in \Q} \{ 1-\phi \geq \varepsilon \}
=\{1-\phi\neq 0\}, \] \noindent{}and by our assumption
($\underline{\lim }\max_{s \in X_{m}}\xi_{m}(s)<1$) we choose $\phi$ such that
$\{1-\phi(z) = 0\}$ is a null set with respect to
the Lebsegue measure. This complete the proof of theorem 3.1. when
$\underline{\lim }\max_{s \in
  X_{m}}\xi_{m}(s)<1.$
$\hfill \Box$\\

\section{On the Ornstein probability space for which $\underline{\lim }\max_{s \in
  X_{m}}\xi_{m}(s)=1$}

Using the same ideas as in the previous section, we have the
following

\begin{lem}
$\displaystyle \limsup_{m \longrightarrow \infty} \iint ||P_m|^2-1|d\lambda d\Prob \geq \iint Q d\lambda d\Prob .$
\end{lem}

\begin{proof}
\noindent{}We have
\begin{eqnarray}\label{eqn:cle}
\iint Q \left | \left |P_m \right|^2-1 \right |d\lambda d\Prob
\geq \int \!\!\left |\int Q \left ( |P_m|^2-1\right )d\Prob \right
| d\lambda,
\end{eqnarray}
\noindent But, from $(4.2)$
\begin{eqnarray}\label{eqn:cle2}
\int \!\!\left (|P_m|^2-1 \right ) d\Prob=2 \re {\left\{\left ( G_{p_m}(z^{h_m+t_m})\right)\left(\int \tau_1 d\Prob
\right) \right \}}\\+ \nonumber
\left | F_{p_m}(z^{h_m+t_m})-\frac{p_m-1}{p_m} \right | \phi_m(z).
\end{eqnarray}
\noindent{}Where, $F_p$ and $G_p$ is define, for any $p \in \N^*$, by
\begin{eqnarray*}
F_p(z)&=& \left |\frac1{\sqrt{p}}\sum_{k=1}^{p-1}z^k \right|^2,\\
G_p(z)&=&\frac1{p}\sum_{k=1}^{p-1}z^k.
\end{eqnarray*}
\noindent{}{\re($z$)} is a real part of the complex number $z$. Combined $(\ref{eqn:cle})$ and $(\ref{eqn:cle2})$ to obtain
\begin{eqnarray}\label{eqn:cle3}
\nonumber \iint Q \left |\left |P_m \right|^2-1 \right |d\lambda
d\Prob & \geq & \iint Q \left ( \left |
F_{p_m}(z^{h_m+t_m})-\frac{p_m-1}{p_m} \right | \phi_m(z)
\right ) d\lambda \Prob -\\
& & 2 \iint Q  \left | G_{p_m}(z^{h_m+t_m})\right|\left|\int \tau_1 d\Prob
\right| d\lambda d\Prob.
\end{eqnarray}
\noindent{}But, on one hand, we have
\begin{eqnarray*}
\int Q  \left | G_{p_m}(z^{h_m+t_m})\right|\left|\int \tau_1 d\Prob
\right| d\lambda &\leq&
{\left (\int \left | G_{p_m}(z^{h_m+t_m})\right|^2 d\lambda \right)}^{\frac12}
{\left (\int Q^2 d\lambda \right)}^{\frac12} \\
&\leq& {\left (\int \left | G_{p_m}(z^{h_m+t_m})\right|^2 d\lambda \right)}^{\frac12}
=\frac1{\sqrt{p_m}}\tend{m}{\infty}0.
\end{eqnarray*}
\noindent{}
\noindent{}On the other hand, since $|X_m| \leq t_m$, $\sum_{k \in X_m} \left (\xi_m
\{k\}\right)^2 \tend{m}{\infty}1$ and for any $f \in L^1$, we have
\[
\widehat{f_{(m)}}(n)=\left \{ \begin{array}{ll} 0  & \mbox{ if $n$
is not divisible by $m$}\\
 \widehat{f}\left (\displaystyle \frac{n}{m} \right ) & \mbox{otherwise}
 \end{array}
 \right.
\]
\noindent{} Where $f_{(m)}(z)=f(z^m)$, we get that
$\displaystyle \left |F_{p_m}(z^{h_m+t_m})-\frac{p_m-1}{p_m} \right |\left |\sum_{k \in X_m} \xi_m(k) z^k \right|^2
d\lambda$ converge to $K. \lambda$, with $K \geq 1$. In fact
\begin{eqnarray*}
&&\int \left |F_{p_m}(z^{h_m+t_m})-\frac{p_m-1}{p_m} \right |\left |\sum_{k \in X_m} \xi_m(k) z^k \right |^2
d\lambda \\ & = &
\sum_{k \in X_m} \left (\xi_m
\{k\}\right)^2
\int \left |F_{p_m}(z^{h_m+t_m})-\frac{p_m-1}{p_m} \right |d\lambda \\
& \geq & \sum_{k \in X_m} \left (\xi_m
\{k\}\right)^2   \int
\left (F_{p_m}(z^{h_m+t_m})-\frac{p_m-1}{p_m} \right ) z^{h_m+t_m} d\lambda \\
 &=&\sum_{k \in X_m}\left (\xi_m
\{k\}\right)^2   \left (\frac{p_m-2}{p_m} \right) \tend{m}{\infty} 1.
\end{eqnarray*}
\noindent{}and the proposition follows from ({\ref {eqn:cle}}).
\end{proof}
\\

\noindent {\bf Proof of the theorem 3.1.in the case of
$\underline{\lim }\max_{s \in
  X_{m}}\xi_{m}(s)=1$}\\

As in the case of $\underline{\lim }\max_{s \in
  X_m}\xi_{m}(s)<1$, we use the lemma (6.1) to etablish  that
  \[
  \lim_{n \longrightarrow \infty} \int \prod_{k=1}^{n}|P_k(z)|
  d\lambda d\Prob=0.
  \]
\noindent{} and the proof of the theorem 3.1. is complete.\\

\begin{rem}

 We note  that Rudolph construction in \cite{Rudolph} is
strictly included in the theory of generalized random Ornstein
construction.
\end{rem}

\noindent{}{\centerline {\bf Acknowledgements}}
 The authors would like to express thanks to J-P. Thouvenot who posed them the
problem of singularity of the spectrum of the Generalized Ornstein
transformations.

\bibliographystyle{nyjplain}

\end{document}